\documentclass[11pt]{article}
\usepackage{}
\usepackage{graphicx}
\usepackage{graphicx,subfigure}
\usepackage{epsfig}
\usepackage{amssymb,pict2e,color}
\usepackage{mathrsfs}
\usepackage{amsmath}
\usepackage{ifpdf}

\newcommand{\E}{{\cal E}}

\newtheorem{theorem}{Theorem}[section]

\newtheorem{lemma}[theorem]{Lemma}
\newtheorem{corollary}[theorem]{Corollary}

\newtheorem{question}[theorem]{Question}

\def\whitebox{{\hbox{\hskip 1pt
 \vrule height 6pt depth 1.5pt
 \lower 1.5pt\vbox to 7.5pt{\hrule width
    3.2pt\vfill\hrule width 3.2pt}%
 \vrule height 6pt depth 1.5pt
 \hskip 1pt } }}
\def\qed{\ifhmode\allowbreak\else\nobreak\fi\hfill\quad\nobreak
     \whitebox\medbreak}

\newcommand{\ignore}[1]{}

\begin {document}

\baselineskip 16pt
\title{On the sharp lower bounds of Zagreb indices of graphs with given number of cut vertices}

 \author{\small   Shengjin Ji$^{a,c}$,  Shaohui Wang$^{b}$\thanks{Corresponding author. \newline
  \E-mail:  jishengjin2013@163.com(S. Ji), shaohuiwang@yahoo.com(S. Wang).}\\
 \small  $^{a}$ School of Science, Shandong University of Technology,
 Zibo, Shandong 255049, China\\
 \small  $^{b}$ Department of Mathematics, Savannah State University, Savannah, GA 31404, USA\\
 \small $^{c}$ School of Mathematics, Shandong University,
Jinan, Shandong 250100, China\\}

\date{}
\maketitle
\begin{abstract}

The first Zagreb index of  a graph $G$ is the sum of the square of every vertex  degree, while the second Zagreb index is the sum of the product of vertex degrees of each edge over all edges.   In our work,  we solve an open question about  Zagreb indices of   graphs  with given number of cut vertices.  The sharp lower bounds are obtained for these indices of graphs in $\mathbb{V}_{n,k}$, where $\mathbb{V}_{n, k}$ denotes the set of all  $n$-vertex graphs with  $k$ cut vertices and at least one cycle.  As consequences, those  graphs with the   smallest   Zagreb indices are characterized.

\end{abstract}

\section{ Introduction}
A topological index is a constant which can be   describing some properties of  a molecular graph,  that is,  a finite  graph  represents the carbon-atom skeleton
of an organic molecule of a hydrocarbon.
During past few decades these have been used for the study of quantitative structure-property relationships (QSPR) and quantitative
structure-activity relationships (QSAR)  and for  the structural essence of biological and chemical compounds.

One of the most well-known topological indices is called Randi\'{c} index,   a  moleculor quantity of branching index \cite{1}. It is known as the classical Randi\'{c} connectivity index, which is the most useful structural descriptor in QSPR and QSAR, see \cite{2,3,4,5}.  Many mathematicians focus considerable interests in the structural and applied issues of  Randi\'{c} connectivity index, see \cite{6,7,8,9}. Based on these perfect considerations,  Zagreb indices\cite{10} are introduced as  expressing formulas for the total $\pi$-electron energy of conjugated molecules below.
\begin{eqnarray} \nonumber
M_1(G) = \sum_{u \in V(G)} d(u)^2
~\text{ and}
~ M_2(G) = \sum_{uv \in E(G)} d(u)d(v),
\end{eqnarray}
where $G$ is a (molecular) graph, $uv$ is a bond between two atoms $u$ and $v$, and $d(u)$ (or $d(v)$, respectively) is the number of atoms that are connected with $u$ (or $v$, respectively).



Zagreb indices are employing as molecular descriptors in QSPR and QSAR, see \cite{11,12}.   In the interdisplinary of mathemactics, chemistry and physics, it is not surprising that  there are numerous studies of properties of the  Zagreb indices of molecular graphs \cite{13,14,1000,1001,1003,1004,01004,WW,1006}.
In \cite{16,17}, some bounds of (chemical) trees on Zagreb indices are studied and surveyed.
Hou et al. \cite{18} found sharp bounds for Zagreb indices of maximal outerplanar graphs.
Li and Zhou \cite{1002} investigated  the maximum and minimum Zagreb indices of graphs with connectivity at most $k$.
 The upper
bounds on Zagreb indices of trees in terms of domination number is studied by Borovi\'canin et al.~\cite{Borov2016}.   In many mathematical literatures~\cite{Bojana2015},   the maximum and minimum
Zagreb indices of trees with a given
number of vertices of maximum degree are explored. Xu and Hua \cite{Xu2012} provided a unified approach to extremal
multiplicative Zagreb indices for trees, unicyclic and bicyclic
graphs.
 Sharp upper and lower bounds of these indices about
$k$-trees are introduced by Wang and Wei~\cite{Wang2015}.
Liu and Zhang  provided several sharp upper bounds for
multiplicative Zagreb indices in terms of graph parameters
such as the order, size and radius~\cite{Liuz20102}.
  The bounds for
the moments and the probability generating function of multiplicative Zagreb indices in a randomly chosen molecular graph with tree structure.
 Zhao and Li \cite{zhaoli} investigated the upper bounds of Zagreb indices, and proposed an open question:

\begin{question}\cite{zhaoli}
How can we determine the lower bound for the first and the second Zagreb
indices of n-vertex connected graphs with k cut vertices? What is the characterization of the
corresponding extremal graphs?
\end{question}

In the view of  above results and open problem,  we proceed to
investigate   Zagreb indices of graphs with  given number of cut vertices in this paper.  It is known that there are many results about Zagreb indices on the graph without cycles. We consider  the set of all  $n$-vertex graphs with  $k$ cut vertices and at least one cycle, denoted by $\mathbb{V}_{n, k}$.  In addition, the  minimum values of $M_1(G) $ and $M_2(G) $ of graphs with   given number of cut vertices are provided.  Furthermore, we characterize  graphs with the  smallest  Zagreb indices in $\mathbb{V}_{n,k}$.

\section{Preliminary}
In this section, we provide some important statements, and introduce several graph transformations. These are significant in the following section.

Let  $G = (V, E)$ be  a simple connected graph of $n$ vertices and $m$ edges, where $V = V (G)$ is  a vertex set and $E = E(G)$ is an edge set. If  $v \in V(G)$, then $N(v)$ is the neighborhood of $v$, that is,  $N_{G}(v)= \left\{ u|\; uv\in E(G)\right\}$  and   the degree of $v$ is $d_{G}(v)= \left| N_G(v)\right|$.
A pendent vertex is the vertex of degree $1$   and a supporting vertex is the vertex adjacent to at least $1$ pendent vertex. A pendent edge is incident to    a supporting vertex and a pendent vertex.
Given sets $S \subseteq V(G)$ and  $F \subseteq E(G)$,   denote by $G[S]$  the subgraph of $G$ induced by $S$, $G - S$   the subgraph induced by $V(G) - S$ and $G - F$   the subgraph of G obtained by deleting $F$.  A vertex $u$ (or an edge $e$, respectively) is said to be a cut vetex (or cut edge, respectively) of a connected graph $G$, if $G- v$ (or $G -e$) has at least two components.   A graph G is called $2$-connected  if there does not exist a vertex whose removal disconnects the graph. A block is a connected graph which does not have any cut vertex. In particular,  $K_2$ is a trivial block, and the endblock contains at most one cut vertex.
Let $P_n$, $S_n$ and $C_n$ be a path, a star and a cycle on $n$ vertices, respectively. Let $T$ be a tree, and $C_m$ be a cycle of $G$. If $G$ contains $T$ as its subgraph via attaching some vertex of $T$ to some vertex of $C_m$, then we say tree $T$ is a pendent tree of $G$. Especially, replacing $T$ by $P_{|T|}$, and choosing its pendent vertex to attach some vertex of $C_m$, we call path $P_{|T|}$ is a pendent path of $G$.  In this exposition we may use the notations and terminology
 of (chemical) graph theory (see \cite{BB,NT}).

We start with an elementary lemma below.

\begin{lemma}\label{l1}
Let $G$ be a graph. If $uv \in E(G)$, then $M_i(G-uv) < M_i(G)$ with $i = 1, 2$.
\end{lemma}

Besides the above lemma, we provide an useful tool about maximal $2$-connected block on Zagreb indices.

\begin{lemma}\label{cyclelemma} Let $G\in \mathbb{V}_n^k$ be a graph with the smallest Zagreb indices  and $D$ a maximal $2$-connected block of $G$ with $i=1, 2$.
If $|D| \geq 3$, then $D$ is a cycle.
\end{lemma}

\begin{proof} If $|D| = 3$, then $D$ is a cycle. Otherwise,
we prove the case of $|D| \geq 4$ by a contradiction. Assume that $D$ is a connected graph without cut vertices and $D$ is not a cycle. Then there exists an edge $uv$ in $D$  such that $D- uv$ has no cut vertices. Obviously, $G- uv \in \mathbb{V}_n^k$.  By Lemma \ref{l1}, $M_i(G-uv) < M_i(G)$, which is contradicted to the choice of $G$.\qed
\end{proof}

The four crucial operations on graphs are given as follows.

\noindent{\bf Operation I.}  As shown in Fig.1, let $H_1$ be a connected graph with $d_{H_1}(v)\geq 3$ and $d_{H_1}(v_1)=1$, and $u_1u_2$ belong to a cycle of $H_1.$
If $H_2=H_1-\{u_1u_2,v_1v\}+\{u_1v_1, u_2v_1\},$ we say that $H_2$ is obtained from $H_1$ by \emph{Operation I}.
\begin{center}
\begin{picture}(322.6,61.9)\linethickness{0.8pt}
\Line(263.7,36.4)(287.7,40.4)
\Line(263.7,36.4)(286.2,29.4)
\put(287.2,35.9){\circle{30.1}}
\put(197.5,12.5){$u_2$}
\put(184.5,44){$u_1$}
\put(204.4,28){$v_1$}
\put(258.9,27.9){$v$}
\put(290.9,40.5){$v_2$}
\put(288.9,25){$v_{\ell}$}
\put(1.5,41.5){$u_1$}
\put(3.5,17){$u_2$}
\put(32.2,34.9){\circle{39}}
\put(71.4,26.4){$v$}
\Line(76.2,34.9)(98.7,27.9)
\Line(76.2,34.9)(83.2,49.4)
\Line(76.2,34.9)(100.2,38.9)
\put(83.4,51.5){$v_1$}
\put(99.7,34.4){\circle{30.1}}
\put(101.4,23.5){$v_{\ell}$}
\put(103.4,39){$v_2$}
\put(218.2,37.4){\circle{42}}
\color[rgb]{0.75,0.75,0.75}\polygon*(239.7,36.9)(264.2,36.9)(264.2,54.9)(239.7,54.9)\color{black}\polygon(239.7,36.9)(264.2,36.9)(264.2,54.9)(239.7,54.9)
\color[rgb]{0.75,0.75,0.75}\polygon*(52.2,34.9)(76.7,34.9)(76.7,52.9)(52.2,52.9)\color{black}\polygon(52.2,34.9)(76.7,34.9)(76.7,52.9)(52.2,52.9)
\put(52,8){$H_1$}
\put(242,8){$H_2$}
\put(239.2,36.4){\circle*{4}}
\put(263.7,36.4){\circle*{4}}
\put(287.7,40.4){\circle*{4}}
\put(286.2,29.4){\circle*{4}}
\put(51.7,34.9){\circle*{4}}
\put(13.4,40){\circle*{4}}
\put(17.9,21.7){\circle*{4}}
\put(76.2,34.9){\circle*{4}}
\put(98.7,27.9){\circle*{4}}
\put(83.2,49.4){\circle*{4}}
\put(100.2,38.9){\circle*{4}}
\put(239.2,36.9){\circle*{4}}
\put(197.9,42.7){\circle*{4}}
\put(198.9,29.1){\circle*{4}}
\put(208.6,18.7){\circle*{4}}
\put(160,-10){\makebox(0,0){Fig.1 The graphs using in Operation {I} and Lemma \ref{rmvertex}.}}
\end{picture}
\end{center}

 Based on the above operation, we obtain a lemma below.

  \begin{lemma}\label{rmvertex}
  If $H_2$ is obtained from $H_1$ by Operation I as shown in Fig.1. Then $M_{i}(H_2)<M_{i}(H_1)$ for $i=1,2.$
  \end{lemma}
  \begin{proof}
  Let $v$ be a vertex of $H_1$ with $d_{H_1}(v)\geq 3$ and containing at least one pendent vertex $v_1$, and $u_1u_2$ be an edge of some cycle in $H_1$ with $d_{H_1}(u_1),d_{H_1}(u_2)\geq 2.$  The neighbors of $v$ are marked as $v_1,v_2,\ldots,v_{\ell}$ for $\ell\geq 3$ (see Fig.1).

   If $v$ doesn't belong to any cycle of $H_1$. Then $H_2$ denotes the graph obtained from $H_1$ by deleting two edges $vv_1,u_1u_2$ and adding edges $u_1v_1,u_2v_1.$ Note that the function $f(x,y)\triangleq xy-x-y+3,$ for $(x,y)\in [2,+\infty)\times [2,+\infty),$ is more than zero. We now deduce that
  \begin{equation*}
  \begin{split}
  M_1(H_1)-M_1(H_2)&=(d_{H_1}(v))^2+(d_{H_1}(v_1))^2
  -(d_{H_2}(v))^2-(d_{H_2}(v_1))^2\\
  &=(d_{H_1}(v)+d_{H_2}(v))-(d_{H_1}(v_1)+d_{H_2}(v_1))\\
  &\geq 5-3=2>0.
  \end{split}
  \end{equation*}
  In terms of the property of $f(x,y),$ for $M_2$, we arrive at
  \begin{equation*}
  \begin{split}
  &M_2(H_1)-M_2(H_2)\\&=\sum\limits_{j=1}^{\ell}d_{H_1}(v)d_{H_1}(v_j)
  +d_{H_3}(u_)d_{H_1}(u_2)\\
  &\ -\sum\limits_{j=2}^{\ell}d_{H_2}(v)d_{H_2}(v_j)
  -d_{H_2}(u_1)d_{H_2}(v_1)-d_{H_2}(u_2)d_{H_2}(v_1)\\
  &=\sum\limits_{j=2}^{\ell}d_{H_1}(v_j)+d_{H_1}(u_1)d_{H_1}(u_2)
  +d_{H_1}(v)-d_{H_2}(u_1)-d_{H_2}(u_2)\\
  &>d_{H_1}(u_1)d_{H_1}(u_2)-d_{H_1}(u_1)-d_{H_1}(u_2)+3\\
  &=f(d_{H_1}(u_1),d_{H_1}(u_2))>0.
  \end{split}
  \end{equation*}
  The special case $v$ belongs to some cycles of $H_1$ should be discussed.
   If $v_1$ is the unique pendent vertex of $H_1$. Then there are nothing to do. If $H_1$ has another pendent vertex, marked as $w_1$, and $H_2=H_1-vv_1+v_1w_1.$ Then the conclusion is also verified. The proof precess of the case is similar with the above argument, so it is omitted.

   Hence, the proof is finished.\qed
  \end{proof}

\noindent{\bf Operation I{I}.} As shown in Fig. 2, let $H_3$ be a graph with $d_{H_3}(v)\geq 3,$ and $w_1w_2$ be an edge included in some cycle of $H_3.$  If $H_4=H_3-\{vv_2,u_{21}v_2,w_1w_2\}+\{v_2w_1,v_2w_2,u_{21}u_{\ell t_{\ell}}\}$ for some $\ell$, we say that $H_4$ is obtained from $H_3$ by \emph{Operation I{I}}.
\begin{center}
\begin{picture}(381.5,72.6)\linethickness{0.8pt}
\put(224.9,34.4){\circle{40}}
\Line(244.9,33.9)(259.9,33.9)
\Line(259.9,33.9)(278.4,33.9)
\Line(278.4,33.9)(278.4,33.9)
\Line(278.4,33.9)(312.9,55.9)
\Line(278.4,33.9)(313.9,34.4)
\put(213.5,18){\small$w_2$}
\put(210,43.5){\small$w_1$}
\put(255.1,25){\small$v_1$}
\put(209,28){\small$v_2$}
\put(290,25){\small$v_{\ell}$}
\put(308,25){\small$u_{\ell t_{\ell}}$}
\put(318.6,38.5){\small$u_{21}$}
\put(341.1,38.5){\small$u_{2t_2}$}
\put(278,25){\small$v$}
\put(23,33){\small$w_1$}
\put(25,18){\small$w_2$}
\put(38.4,32.4){\circle{40}}
\Line(58.4,31.9)(73.4,31.9)
\put(68.6,25){\small$v_1$}
\Line(73.4,31.9)(91.9,31.9)
\put(99,53){\small$u_{21}$}
\put(91,42){\small$v_2$}
\put(92.1,25){\small$v$}
\put(105,64){\small$u_{2t_2}$}
\Line(91.9,31.9)(142.4,32.4)
\Line(91.9,31.9)(129.4,48.4)
\Line(91.9,31.9)(122.9,60.4)
\put(105.1,25){\small$v_{\ell}$}
\put(122.1,25){\small$u_{\ell 1}$}
\put(138.1,25){\small$u_{\ell t_{\ell}}$}
\Line(313.9,34.4)(330.9,34.9)
\Line(330.9,34.9)(352.9,34.9)
\put(68,1){$H_3$}
\put(256.5,0){$H_4$}
\polygon(200.9,9.4)(248.9,9.4)(248.9,58.4)(200.9,58.4)
\polygon(11.9,9.9)(62.9,9.9)(62.9,56.9)(11.9,56.9)
\put(244.9,33.9){\circle*{4}}
\put(259.9,33.9){\circle*{4}}
\put(278.4,33.9){\circle*{4}}
\put(312.9,55.9){\circle*{4}}
\put(313.9,34.4){\circle*{4}}
\put(290.4,41.6){\circle*{4}}
\put(301.6,48.7){\circle*{4}}
\put(292.1,34.1){\circle*{4}}
\put(209.5,47.3){\circle*{4}}
\put(216.8,16.1){\circle*{4}}
\put(18.7,36.3){\circle*{4}}
\put(58.4,31.9){\circle*{4}}
\put(22.5,20.2){\circle*{4}}
\put(73.4,31.9){\circle*{4}}
\put(91.9,31.9){\circle*{4}}
\put(142.4,32.4){\circle*{4}}
\put(129.4,48.4){\circle*{4}}
\put(122.9,60.4){\circle*{4}}
\put(103,42.2){\circle*{4}}
\put(105,37.7){\circle*{4}}
\put(111.3,32.1){\circle*{4}}
\put(112.5,50.9){\circle*{4}}
\put(112.5,50.9){\circle*{4}}
\put(117.1,43){\circle*{4}}
\put(127.3,32.3){\circle*{4}}
\put(205.3,30.4){\circle*{4}}
\put(330.9,34.9){\circle*{4}}
\put(352.9,34.9){\circle*{4}}
\put(180,-15){\makebox(0,0){Fig.2 The graphs using in Operation I{I} and Lemma \ref{rmspath}.}}
\end{picture}
\end{center}

 \begin{lemma}\label{rmspath}
  If $H_4$ is obtained from $H_3$ by Operation I{I} as shown in Fig.2. Then $M_{i}(H_4)<M_{i}(H_3)$ for $i=1,2.$
  \end{lemma}
  \begin{proof}
  Let $v\in V(H_3)$ with $d_{H_3}(v)\geq3$ and $w_1w_2$ be an edge of some cycle in $H_3$. Its neighbors are labeled as $v_1,v_2,\ldots,v_{\ell}(\ell\geq 3).$ If there is at least one pendent vertex of $v$, then this case can be reduced to Operation I. So we may assume that $v$ only possesses pendent paths whose length are no less than $2$. If $v$ has an unique pendent path, there are nothing to do. If $v$ contains at least two pendent paths, e.g., $P_2(v_2 u_{21} \ldots u_{2t_2})$ and $P_\ell(v_\ell \ldots u_{\ell 1}u_{\ell t_{\ell}})$ with $t_2,t_{\ell}\geq 1$. Let $H_4=H_3-\{vv_2,u_{21}v_2,w_1w_2\}+\{v_2w_1,v_2w_2,u_{21}u_{\ell t_{\ell}}\}.$ Observe that the function $g(x,y)\triangleq xy-2x-2y+5,$ for $(x,y)\in [2,+\infty)\times [2,+\infty),$ is more than zero.
  We now deduce that
  \begin{equation*}
  \begin{split}
  M_1(H_3)-M_1(H_4)&=(d_{H_3}(v))^2+(d_{H_3}(u_{\ell}t_{\ell}))^2
  -(d_{H_4}(v))^2-(d_{H_4}(u_{\ell}t_{\ell}))^2\\
  &=(d_{H_3}(v)+d_{H_4}(v))
  -(d_{H_3}(u_{\ell}t_{\ell})+d_{H_4}(u_{\ell}t_{\ell}))\\
  &\geq 5-3=2>0.
  \end{split}
  \end{equation*}
  Since $d_{H_3}(w_1),d_{H_3}(w_2)\geq 2.$ For $M_2$, in terms of the property of $g(x,y),$ we get
  \begin{equation*}
  \begin{split}
  &M_2(H_3)-M_2(H_4)\\&=\sum\limits_{j=1}^{\ell}d_{H_3}(v)d_{H_3}(v_j)
  +d_{H_3}(w_1)d_{H_3}(w_2)+d_{H_3}(u_{21})d_{H_3}(v_2)
  +d_{H_3}(u_{\ell (t_{\ell}-1)})d_{H_3}(u_{\ell t_{\ell}})\\
  &\quad-\sum\limits_{j=1,j\neq 2}^{\ell}d_{H_4}(v)d_{H_4}(v_j)
  -d_{H_4}(v_2)d_{H_4}(w_1)-d_{H_4}(v_2)d_{H_4}(w_2)\\
  &\quad-d_{H_4}(u_{\ell (t_{\ell}-1)})d_{H_4}(u_{\ell t_{\ell}})-d_{H_4}(u_{21 })d_{H_4}(u_{\ell t_{\ell}})\\
  &=\sum\limits_{j=1,j\neq 2}^{\ell}d_{H_3}(v_j)+d_{H_3}(w_1)d_{H_3}(w_2)+2d_{H_3}(v)
  -2d_{H_3}(w_1)-2d_{H_3}(w_2)-2\\
  &\geq d_{H_3}(w_1)d_{H_3}(w_2)
  -2d_{H_3}(w_1)-2d_{H_3}(w_2)+6\\
  &=g(d_{H_3}(w_1),d_{H_3}(w_2))+1>0.
  \end{split}
  \end{equation*}Hence, the conclusion is verified. \qed
  \end{proof}

 \noindent{\bf Operation I{I}{I}.} As shown in Fig. 3, let $G_0$ be a connected graph with $|G_0|\geq 2$ and having two vertices $u$ and $w$, and $G_1$ be the graph which contains a cycle $C_1.$  Let $H_5$ be a graph on order $n(\geq 6)$ obtained from $G_0$ by identifying some vertex of $C_1$ with vertex $u$ and some vertex of $C_2$ with vertex $w$, respectively.
If $H_6$ denote the new graph from $H_5-\{v_1w,v_0v_2,u_1u_2\}+\{u_1v_0,u_2v_1\},$ we say that $H_6$ is obtained from $H_5$ by \emph{Operation I{I}{I}}.

\begin{center}
\begin{picture}(297.6,68.4)\linethickness{0.8pt}
\put(216.7,38.9){\circle{49}}
\cbezier(241.2,40.4)(242,19.9)(264.9,19.9)(265.7,40.4)\cbezier(241.2,40.4)(242,61)(264.9,61)(265.7,40.4)
\put(232,35.9){\small$u$}
\put(266,32){\small$w$}
\put(190,18){\small$v_1$}
\put(243.9,43){\small$v_3$}
\put(190,56){\small$v_0$}
\put(212.5,54){\small$u_1$}
\put(210,18.5){\small$u_2$}
\put(246,28){\small$v_{\ell}$}
\put(282.4,46){\small$v_2$}
\put(6,46){\small$u_1$}
\put(8,20.5){\small$u_2$}
\put(39.2,38.4){\circle{49}}
\put(54.5,34.9){\small$u$}
\put(68,28){\small$v_{\ell}$}
\put(66.4,42){\small$v_3$}
\cbezier(64.2,39.4)(65,17.5)(88.9,17.5)(89.7,39.4)\cbezier(64.2,39.4)(65,61.3)(88.9,61.3)(89.7,39.4)
\Line(75.2,33.4)(89.7,38.3)
\Line(75.2,46.4)(89.7,38.3)
\put(88,53){\small$v_2$}
\put(89,17){\small$v_1$}
\put(91,33){\small$w$}
\put(108.2,37.4){\circle{37}}
\put(109.4,58){\small$v_0$}
\Line(265.7,39.9)(290.2,39.9)
\Line(252.7,47.4)(265.7,39.9)
\Line(252.7,34.4)(265.7,39.9)
\put(56.8,5){$H_5$}
\put(236.8,5){$H_6$}
\put(197.5,31){$C_1$}
\put(109.5,31.5){$C_2$}
\put(27,31){$C_1$}
\put(241.2,39.4){\circle*{4}}
\put(265.7,39.9){\circle*{4}}
\put(252.7,47.4){\circle*{4}}
\put(252.7,34.4){\circle*{4}}
\put(215.6,63.4){\circle*{4}}
\put(213,14.7){\circle*{4}}
\put(63.7,38.9){\circle*{4}}
\put(15.6,44.9){\circle*{4}}
\put(17.9,26.4){\circle*{4}}
\put(75.2,33.4){\circle*{4}}
\put(89.7,38.3){\circle*{4}}
\put(75.2,46.4){\circle*{4}}
\put(97.1,22.6){\circle*{4}}
\put(97.6,52.6){\circle*{4}}
\put(114.1,55){\circle*{4}}
\put(200.1,20.9){\circle*{4}}
\put(203.5,59.5){\circle*{4}}
\put(290.2,39.9){\circle*{4}}
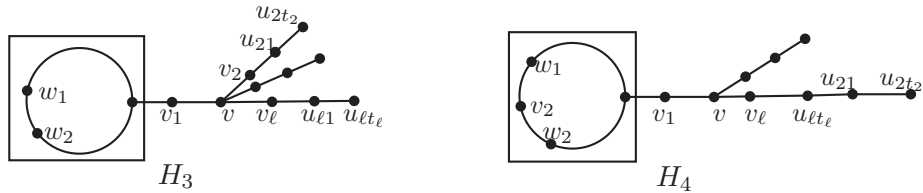
\put(160,-10){\makebox(0,0){Fig.3 The graphs using in Operation I{II} and Lemma \ref{bc1}.}}
\end{picture}
\end{center}

  \begin{lemma}\label{bc1}
  If $H_6$ is obtained from $H_5$ by Operation I{I}{I} as shown in Fig.3. Then $M_{i}(H_6)<M_{i}(H_5)$ for $i=1,2.$
  \end{lemma}
  \begin{proof}
  Let $H_5$ be the graph shown in Fig.3, $u$ and $w$ be two cut-vertex of $H_5.$ $C_1$ and $C_2$ are its two cycles, where $C_2$ is an endblock. $v_1w,v_0v_2,v_2w\in E(C_2)$ and $u_1u_2\in E(C_1)$ with with $d_{H_5}(u_1),d_{H_5}(u_2)\geq 2.$ Let $H_6=H_5-\{v_2v_0,u_1u_2,v_1w\}+\{u_1v_0,u_2v_1\}.$ We can obtain that
  \begin{equation*}
  \begin{split}
  M_2(H_5)-M_2(H_6)&=(d_{H_5}(w))^2+(d_{H_5}(v_2))^2
  -(d_{H_6}(w))^2-(d_{H_6}(v_2))^2\\
  &=(d_{H_5}(w)+d_{H_6}(w))-(d_{H_5}(v_2)+d_{H_6}(v_2))\\
  &\geq 5-3=2>0.
  \end{split}
  \end{equation*}Similarly, For $M_2$, we can deduce that
  \begin{equation*}
  \begin{split}
  &M_2(H_5)-M_2(H_6)\\
  &=\sum\limits_{j=1}^{t}d_{H_5}(w)d_{H_5}(v_j)+d_{H_5}(v_2)d_{H_5}(v_0)
  +d_{H_5}(u_1)d_{H_5}(u_2)\\
  &\quad-\sum\limits_{j=2}^{t}d_{H_6}(w)d_{H_6}(v_j)-d_{H_6}(u_1)d_{H_6}(v_0)
  -d_{H_6}(u_2)d_{H_6}(v_1)\\
  &=\sum\limits_{j=3}^{t}d_{H_5}(v_j)+3d_{H_5}(w)
+d_{H_5}(u_1)d_{H_5}(u_2)-d_{H_5}(u_1)-d_{H_5}(u_2)\\
  &\geq d_{H_5}(u_1)d_{H_5}(u_2)-d_{H_5}(u_1)-d_{H_5}(u_2)+11\\
  &=f(d_{H_5}(u_1),d_{H_5}(u_2))+8>0.
  \end{split}
  \end{equation*}Therefore, the proof is finished.\qed
  \end{proof}
\noindent {\bf Operation {I}V.} As shown in Fig. 4, let $G_0$ be a connected graph having a vertex $v$, and $G_1$ be a graph which contains a cycle $C_1.$ $H_7$ denotes the graph by attaching some vertex of $C_1$ and $C_2$ to the vertex $v$, respectively. Clearly, $C_2$ is an endblock of $H_7.$  If $H_8=H_7-\{vv_2,v_0v_1\}+v_0v_2$, we say that $H_8$ is obtained from $H_7$ by \emph{Operation I{V}}.
\begin{center}
\begin{picture}(214.5,81.4)\linethickness{0.8pt}
\cbezier(40.2,50.4)(3.2,12.4)(73.2,9.4)(40.2,50.9)
\cbezier(39.7,50.9)(91.7,107.4)(88.7,-4.6)(40.2,51.4)
\cbezier(39.2,50.9)(-21.3,6.9)(-2.8,112.4)(39.7,50.9)
\put(30.5,26.5){\scriptsize$G_0$}
\put(46,64){\small$v_1$}
\put(60.4,71){\small$v_0$}
\put(30.9,62.5){\small$v_2$}
\put(43,47.9){\small$v$}
\cbezier(175.7,48.9)(124.2,6.9)(125.2,98.4)(165.7,67.4)
\put(166,24.5){\scriptsize$G_0$}
\cbezier(175.7,48.4)(138.7,10.4)(208.7,7.4)(175.7,48.9)
\put(155.4,73){\small$v_2$}
\cbezier(181.7,66.4)(227.2,104.9)(224.2,-7.1)(175.7,48.9)
\put(179,45.9){\small$v$}
\put(189.4,57){\small$v_1$}
\put(177.9,73){\small$v_0$}
\Line(165.7,67.4)(181.7,66.4)
\Line(190.7,64.4)(175.7,48.9)
\put(28.3,8){$H_7$}
\put(161.8,8){$H_8$}
\put(54.4,45.5){\scriptsize$C_2$}
\put(10.4,49.5){\scriptsize$C_1$}
\Line(39.7,50.9)(36.7,39.4)
\Line(39.7,50.9)(43.2,39.4)
\Line(175.7,48.9)(172.2,37.9)
\Line(175.7,48.9)(178.2,36.9)
\put(40.2,50.4){\circle*{4}}
\put(40.2,50.9){\circle*{4}}
\put(39.7,50.9){\circle*{4}}
\put(40.2,51.4){\circle*{4}}
\put(39.2,50.9){\circle*{4}}
\put(39.7,50.9){\circle*{4}}
\put(31.5,61.1){\circle*{4}}
\put(49.9,60.5){\circle*{4}}
\put(63.6,67.2){\circle*{4}}
\put(175.7,48.9){\circle*{4}}
\put(165.7,67.4){\circle*{4}}
\put(175.7,48.4){\circle*{4}}
\put(175.7,48.9){\circle*{4}}
\put(175.7,48.9){\circle*{4}}
\put(181.7,66.4){\circle*{4}}
\put(190.7,64.4){\circle*{4}}
\put(175.7,48.9){\circle*{4}}
\put(36.7,39.4){\circle*{4}}
\put(43.2,39.4){\circle*{4}}
\put(172.2,37.9){\circle*{4}}
\put(178.2,36.9){\circle*{4}}
\put(110,-10){\makebox(0,0){Fig.4 The graphs using in Operation I{V} and Lemma \ref{bc2}.}}
\end{picture}
\end{center}

  \begin{lemma}\label{bc2}
  If $H_8$ is obtained from $H_7$ by Operation I{V} as shown in Fig.4. Then $M_{i}(H_8)<M_{i}(H_7)$ for $i=1,2.$
  \end{lemma}
  \begin{proof}
As shown in Fig.4, two cycles $C_1$ and $C_2$ of $H_7$ share a common vertex $v$ with $d_{H_7}(v)\geq 4$ whose neighbors are labeled as $v_1,v_2,\ldots v_t.$ Obviously, $t\geq 4.$ In addition, $C_2$ is an endblock of $H_7$. Let $H_8$ denote the new graph obtained from $H_7$ by deleting edges $vv_2,v_0v_1$ and linking $v_2$ to $v_0$. We will deduce the relations of the two graphs $H_7$ and $H_8$ in terms of $M_1$ and $M_2,$ respectively.
\begin{equation}\label{41}
\begin{split}
M_1(H_7)-M_1(H_8)&=(d_{H_7}(v))^2+(d_{H_7}(v_1))^2
  -(d_{H_8}(v))^2-(d_{H_8}(v_1))^2\\
  &=(d_{H_7}(v)+d_{H_8}(v))+(d_{H_7}(v_1)+d_{H_8}(v_1))\\
  &\geq 7+3=10>0,
\end{split}
\end{equation}
\begin{equation}\label{42}
\begin{split}
M_2(H_7)-M_2(H_8)
&=\sum\limits_{j=1}^{t}d_{H_7}(v)d_{H_7}(v_{j})+d_{H_7}(v_0)d_{H_7}(v_1)\\
&\quad -\sum\limits_{j=3}^{t}d_{H_8}(v)d_{H_8}(v_{j})-d_{H_8}(v)d_{H_8}(v_{1})
-d_{H_8}(v_0)d_{H_8}(v_2)\\
&=\sum\limits_{j=3}^{t}d_{H_7}(v_{j})+d_{H_7}(v)+(d_{H_7}(v)-2)d_{H_7}(v_2)
+5>0.
\end{split}
\end{equation}Together Eq.\ref{41} with Eq.\ref{42}, the conclusion is verified. \qed
  \end{proof}

  Let $H$ be a connected graph with $|E(H)|-|V(H)|\geq 0$ and $u,v\in V(H)$ contained in a  cycle of $H$. Denote by $H(a,b)$ the graph formed from $H$ by attaching two paths $P_a$ and $P_b$ to $u$ and $v$, respectively.
 \begin{lemma}\label{rmpath}
  For $d_{H(a,b)}(u),d_{H(a,b)}(v)\geq 3$, we have $M_i(H(a,b))\geq M_i(H(1,a+b-1))$ for $i=1,2.$
 \end{lemma}
 \begin{proof}
 Since $u,v$ belong to some cycle of $H,$ we have $d_{H(a,b)}(u),d_{H(a,b)}(v)\geq 3.$ Without loss of generality, assume that $d_{H(a,b)}(u)\geq d_{H(a,b)}(v).$ We now label all vertices of the two paths $P_a$ and $P_b$
as $uu_1u_2\ldots u_{a-1}$ and $vv_1v_2\ldots v_{b-1},$ respectively. Suppose that, besides $u_1$, the other neighbors of $u$ are $w_1,w_2,\ldots,w_t$ with $t\geq 2.$  $H(1,a+b-1)$ is the graph formed from $H(a,b)$ by deleting edge $uu_1$ and connecting $u_1$ with $v_{b-1}$. For short, we mark $H(a,b)$ and $H(1,a+b-1)$ as $H_0$ and $H'_0,$ respectively.
We first consider $M_1$ and deduce that
\begin{equation*}
\begin{split}
M_1(H_0)-M_1(H'_0)&=(d_{H_0}(u))^2+(d_{H_0})(v_{b-1}))^2
-(d_{H'_0}(u))^2-(d_{H'_0}(v_{b-1}))^2\\
&=d_{H_0}(u)+d_{H'_0}(u)+3>0.
\end{split}
\end{equation*}Similarly, for $M_2$, we get that
\begin{equation*}
\begin{split}
&M_2(H_0)-M_2(H'_0)\\
=&\sum\limits_{j=1}^{t}d_{H_0}(u)d_{H_0}(w_j)
+d_{H_0}(u)d_{H_0}(u_1)+d_{H_0}(v_{b-2})d_{H_0}(v_{b-1})\\
&-\sum\limits_{j=1}^{t}d_{H'_0}(u)d_{H'_0}(w_t)
-d_{H'_0}(v_{b-2})d_{H'_0}(v_{b-1})-d_{H'_0}(v_{b-1})d_{H'_0}(u_1)\\
=&\sum\limits_{j=1}^{t}d_{H_0}(w_j)+d_{H_0}(u_1)(d_{H_0}(u)-d_{H'_0}(v_{b-1}))
-d_{H_0}(v_{b-2})\\
\geq& d_{H_0}(u)+d_{H_0}(u_1)-d_{H_0}(v)>0.
\end{split}
\end{equation*}Therefore, we complete the proof.
\qed
 \end{proof}
  Especially, the two vertices $u$ and $v$ are identified in $H(a,b)$. Then, use the similar way of Lemma \ref{rmpath}, we also got a new graph $H(2,a+b-2)$ such that $M_i(H(a,b))\geq M_i(H(a',b'))$ with $a'=2,b'=a+b-2$ for $i=1,2.$ Obviously, $P_{a'}=uu_1$ and $u_1$ is a pendant. Hence, from Lemma \ref{rmvertex}, we deduce that there exists $H'$ with $|H'|=|H|+1$.(It is obtained from $H$ by subdividing its one edge $w_1w_2$ included in some cycle and marking the vertex as $u_1$.) such that $M_i(H(a,b))\geq M_i(H'(1,a+b-2))$ for $i=1,2.$ We list the result as follows.
  \begin{corollary}\label{rmpaths}
  If  two vertices $u$ and $v$ are identified in $H(a,b)$. Then there exists a graph $H'$ on order $|H|+1$ such that $M_i(H(a,b))\geq M_i(H'(1,a+b-2))$ for $i=1,2.$
  \end{corollary}

\section{Main results}
In this section, we provide the lowest bounds on Zagreb indices of graphs in $\mathbb{V}_n^k$. The corresponding graphs are characterized as well.

\begin{theorem}
Let $G$ be a graph in $\mathbb{V}_n^k$. Then  \\
 (i) $M_1(G)\geq 4n+2$, the equality holds if and only if $G \cong C_{n,k}$,\\
 (ii) $M_2(G)\geq 4n+4$, the equality holds if and only if $G \cong C_{n,k}$.
 \end{theorem}
 \begin{proof}
Choose a graph $G\in \mathbb{V}_n^k$ such that $G$ has the minimal value of $M_i$ with $i=1, 2$.  Let $B$ be a cut vertex set of size $k$ in $G$. $G$ can be divided into $s$ blocks via the $k$ cut vertices, and they are denoted by $D_1, D_2, \cdots, D_s$. Clearly, $|D_j|=2$ or $|D_j|\geq 3$ for some $j$.
 We start with a claim.

 \noindent{\bf Claim 1} $G$ has only one pendent tree. In fact, the tree is a path.

\begin{proof}
Since $G$ is the graph for which $M_i(G)$ has the minimum for $i=1,2$ in all connected graphs possessing $k$ cut vertices. We claim that $G$ includes at least a pendent tree. If not, we will get a new graph $G'$ from $G$, and by Lemma \ref{bc1} Lemma \ref{bc2} and $M_i(G')$ is less than $M_i(G)$. We get a contradiction. In addition, every pendent tree of $G$ must be a path. If not, from Lemma \ref{rmspath}, there exists a new graph $G''$ such that $M_i(G'')<M_i(G)$, which contradicts with the choice of $G.$
If $G$ includes at least two pendent paths. By means of Lemma \ref{rmpath} and Corollary \ref{rmpaths}, there is a graph $G_1$ for which $M_i(G_1)<M_i(G)$. This is a contradiction. Note that the number $|B|$ is not changed during these operations. Thus, we complete the proof of this claim. \qed
\end{proof}

According to Lemma \ref{cyclelemma} and Claim 1,
 we know that $G$ is a block graph and its blocks consists of cycle and $K_2$, and $G$ has a unique pendent path, marked as $X(P)$. If $G$ just contains one cycle, then there is nothing to do. We now suppose that $G$ possesses at least two cycles.

 We now claim that all endblocks of $G$ are cycles except for $K_2$ of $X(P)$. Otherwise, $G$ has no less than two pendent paths which contradicts with Claim 1.

 {\bf Case 1.} $G$ just includes two endblocks.

  According the above argument, we can deduce that the two endblocks of $G$ are one cycle $C_1$ and $K_2.$ From the assumption, $G$ contains another cycle $C_2$.
 In terms of Lemma \ref{bc1} and Lemma \ref{bc2}, there is a graph $G'$ for which $M_i(G')<M_i(G)$ for $i=1,2.$

 {\bf Case 2.} The number of endblocks in $G$ is more than two.

  By means of the assumption, $G$ includes at least two cycles endblocks, e.g., $C_3$ and $C_4$. We will get a new graph $G"$ obtained from $G$
 such that $M_i(G')<M_i(G)$ for $i=1,2$ through Lemma \ref{bc1} and Lemma \ref{bc2}.

 By combining Case 1 and Case 2, we deduce a contradiction with the choice of $G$. Hence, $G$ just possesses unique cycle $C_5$. Since $G$ belongs to $\mathbb{V}_n^k$, we can deduce that $C_5\cong C_{n-k}$ and $X(P)\cong P_{k+2}.$ Therefore, $G\cong C_{n,k}.$ By direct calculation,
 We arrive at $M_1(C_{n,k})=4n+2$,  $M_2(C_{n,k})=4n+4.$ We hence complete the proof. \qed
 \end{proof}


 \vspace{5pt} \noindent
{\bf Acknowledgments}

This work was partial supported by National Natural Science Foundation of China (Grant Nos.11401348 and 11561032), and supported by Postdoctoral Science Foundation of China and the China Scholarship Council.

\end{document}